\documentclass[a4paper,12pt,final]{amsart}
\usepackage{times,a4wide,mathrsfs}

\newcommand{\N}{\mathbb{N}}
\newcommand{\D}{\mathbb{D}}

\newcommand{\T}{\mathbb{T}}

\newtheorem*{theorem}{Theorem}

\begin{document}

\title[A lower bound in Nehari's theorem on the polydisc]{A
lower bound in Nehari's theorem on the polydisc}

\author[J. Ortega-Cerd\`a]{Joaquim Ortega-Cerd\`a}
\address{Dept.\ Matem\`atica Aplicada i An\`alisi,
Universitat  de Barcelona, Gran Via 585, 08071 Barce- lona, Spain}
\email{jortega@ub.edu}

\author{Kristian Seip}
\address{Department of Mathematical Sciences, Norwegian University of Science and Technology,
NO-7491 Trondheim, Norway} \email{seip@math.ntnu.no}
\thanks{The first author is supported by the project
MTM2008-05561-C02-01 and and the grant
2009 SGR 1303. The second author is supported by the Research
Council of Norway grant 160192/V30. This work was done as part of
the research program \emph{Complex Analysis and Spectral Problems}
at Centre de Recerca Matem\`{a}tica (CRM), Bellaterra in the spring
semester of 2011. The authors are grateful to CRM for its support
and hospitality.}



\begin{abstract}
By theorems of Ferguson and Lacey ($d=2$) and Lacey and Terwilleger
($d>2$), Nehari's theorem is known to hold on the polydisc $\D^d$
for $d>1$, i.e., if $H_\psi$ is a bounded Hankel form on $H^2(\D^d)$
with analytic symbol $\psi$, then there is a function $\varphi$ in
$L^\infty(\T^d)$ such that $\psi$ is the Riesz projection of
$\varphi$. A method proposed in Helson's last paper is used to show
that the constant $C_d$ in the estimate $\|\varphi\|_\infty\le C_d
\|H_\psi\|$ grows at least exponentially with $d$; it follows that
there is no analogue of Nehari's theorem on the infinite-dimensional
polydisc.
\end{abstract}

\maketitle

This note solves the following problem studied by H. Helson
\cite{Helbook, H}: Is there an analogue of Nehari's theorem on the
infinite-dimensional polydisc? By using a method proposed in
\cite{H}, we show that the answer is negative. The proof is of
interest also in the finite-dimensional situation because it gives a
nontrivial lower bound for the constant appearing in the norm
estimate in Nehari's theorem; we choose to present this bound as our
main result.

We first introduce some notation and give a brief account of
Nehari's theorem. Let $d$ be a positive integer, $\D$ the open unit
disc, and $\T$ the unit circle. We let $H^2(\D^d)$ be the Hilbert
space of functions analytic in $\D^d$ with square-summable Taylor
coefficients. Alternatively, we may view $H^2(\D^d)$ as a subspace
of $L^2(\T^d)$ and express the inner product of $H^2(\D^d)$ as
$\langle f,g\rangle=\int_{\T^d} f\overline{g}$, where we integrate
with respect to normalized Lebesgue measure on $\T^d$. Every
function $\psi$ in $H^2(\D^d)$ defines a Hankel form $H_\psi$ by the
relation $H_\psi(fg)=\langle fg, \psi\rangle$; this makes sense at
least for holomorphic polynomials $f$ and $g$. Nehari's theorem---a
classical result \cite{Ne} when $d=1$ and a remarkable and
relatively recent achievement of S. Ferguson and M. Lacey \cite{FL}
($d=2$) and M. Lacey and E. Terwilleger \cite{LaTe} ($d>2$) in the
general case---says that $H_\psi$ extends to a bounded form on
$H^2(\D^d)\times H^2(\D^d)$ if and only if $\psi=P_+\varphi$ for
some bounded function $\varphi$ on $\T^d$; here $P_+$ is the Riesz
projection on $\T^d$ or, in other words, the orthogonal projection
of $L^2(\T^d)$ onto $H^2(\D^d)$. We define $C_d$ as the smallest
constant $C$ that can be chosen in the estimate
\[ \|\varphi\|_\infty\le C \| H_\psi\|,\]
where it is assumed that $\varphi$ has minimal $L^\infty$ norm.
Nehari's original theorem says that $C_1=1$.

\begin{theorem}
For even integers $d\ge 2$, the constant $C_{d}$ is at least
$(\pi^2/8)^{d/4}$.
\end{theorem}

The theorem thus shows that the blow-up of the constants observed in
\cite{La, LaTe} is not an artifact resulting from the particular
inductive argument used there.

Since clearly $C_d$ increases with $d$ and, in particular, we would
need that $C_d\le C_\infty$ should Nehari's theorem extend to the
infinite-dimensional polydisc, our theorem gives a negative solution
to Helson's problem.

Nehari's theorem can be rephrased as saying that functions in
$H^1(\D^d)$ (the subspace of holomorphic functions in $L^1(\T^d)$)
admit weak factorizations, i.e., every $f$ in $H^1(\D^d)$ can be
written as $f=\sum_j g_j h_j$ with $f_j$, $g_j$ in $H^2(\D^d)$ and $
\sum_j \|g_j\|_2 \|h_j\|_2\le A \|f\|_1$ for some constant $A$.
Taking the infimum of the latter sum when $g_j$, $h_j$ vary over all
weak factorizations of $f$, we get an alternate norm (a projective
tensor product norm) on $H^1(\D^d)$ for which we write
$\|f\|_{1,w}$. We let $A_d$ denote the smallest constant $A$ allowed
in the norm estimate $\|f\|_{1,w}\le A \|f\|_1$. Our proof shows
that we also have $A_d\ge (\pi^2/8)^{d/2}$ when $d$ is an even
integer.

\begin{proof}[Proof of the theorem]

We will follow Helson's approach \cite{H} and also use his
multiplicative notation. Thus we define a Hankel form on $\T^\infty$
as
\[ H_{\psi}(fg)=\sum_{j,k=1}^\infty \rho_{jk} a_j b_k; \]
here $(a_j)$, $(b_j)$, and $(\rho_j)$ are the sequences of
coefficients of the power series of the functions $f$, $g$, and
$\psi$, respectively. More precisely, we let $p_1$, $p_2$, $p_3$,
... denote the prime numbers; if $j=p_1^{\nu_1}\cdots p_k^{\nu_k}$,
then $a_j$ (respectively $b_j$ and $\rho_j$) is the coefficient of
$f$ (respectively of $g$ and $\psi$) with respect to the monomial
$z_1^{\nu_1}\cdots z_k^{\nu_k}$. We will only consider the
finite-dimensional case, which means that the coefficients will be
nonzero only for indices $j$ of the form $p_1^{\nu_1}\cdots
p_d^{\nu_d}$. The prime numbers will play no role in the proof
except serving as a convenient tool for bookkeeping.

We now assume that $d$ is an even integer and introduce the set
\[ I=\left\{n\in\N:\ n=\prod_{j=1}^{d/2} q_j \ \ \text{and} \ \ q_j=p_{2j-1} \ \ \text{or}\ \
q_j=p_{2j}\right\}.\]  We define a Hankel form $H_\psi$ on $\D^d$ by
setting $\rho_n=1$ if $n$ is in $I$ and $\rho_n=0$ otherwise.

We follow \cite[pp. 81--82]{H} and use the Schur test to estimate
the norm of $H_\psi$. It suffices to choose a suitable finite
sequence of positive numbers $c_j$ with $j$ ranging over those
positive integers that divide some number in $I$; for such $j$ we
choose
\[ c_j=2^{-\Omega(j)/2},\]
where $\Omega(j)$ is the number of prime factors in $j$. We then get
\[ \sum_{k} \rho_{jk} c_k=2^{d/2-\Omega(j)} \cdot 2^{-(d/2-\Omega(j))/2}=2^{d/4} c_j,\]
so that $\|H_\psi\|\le 2^{d/4}$ by the Schur test.

If $f$ is a function in $H^1(\D^d)$ with associated Taylor
coefficients $a_n$, then
\[
H_\psi (f)=  \sum_n a_n \rho_n. \] We choose
\begin{equation}\label{ourf} f(z)=\prod_{j=1}^{d/2} (z_{2j-1}+z_{2j})
\end{equation}
for which $a_n=\rho_n$ and thus $H_\psi(f)=2^{d/2}$. On the other
hand, an explicit computation shows that
\[ \|f\|_1=(4/\pi)^{d/2} \]
so that $H_\psi$, viewed as a linear functional on $H^1(\D^d)$, has
norm at least $(\pi/2)^{d/2}$. This concludes the proof since it
follows that we must have $(\pi/2)^{d/2}\le \|\varphi\|_\infty$ and
we know from above that $\|H_\psi\|\le 2^{d/4}$.
\end{proof}
It is worth noting that our application of the Schur test shows that
in fact $\|H_\psi\|=2^{d/4}$ since $\| f\|_2=2^{d/4}$. The fact that
$|H_\psi (f)|=\|H_\psi\| \|f\|_2$ implies that
\[ \| f\|_{1,w}=\|f\|_2.\]
In other words, the trivial factorization $f\cdot 1$ is an optimal
weak factorization of the function $f$ defined in \eqref{ourf}.

\end{document}